\documentclass[12pt]{article}
\usepackage{natbib}
\usepackage{hyperref}
\usepackage{amssymb}
\usepackage{dsfont}
\usepackage{hhline}
\begin{document}
\newtheorem{theorem}{Theorem}
\newtheorem{remarque}{Remark}
\newtheorem{lemma}{Lemma}
\newtheorem{corollaire}{Corollary}

%
%
\vskip 3mm

\vskip 3mm \baselineskip=15pt

\noindent ADDITIVE REGRESSION MODEL FOR CONTINUOUS TIME PROCESSES

\vskip 5mm \noindent Mohammed DEBBARH and Bertrand MAILLOT \vskip
2mm \noindent

\noindent Universit\'e Paris 6

\noindent 175,\ Rue du Chevaleret,  75013 Paris.

\noindent debbarh@ccr.jussieu.fr.
\vskip 4mm \baselineskip=18pt

\noindent ABSTRACT

In the setting of additive regression model for continuous time
process, we establish the optimal uniform convergence rates and
optimal asymptotic quadratic error of additive regression. To
build our estimate, we use the marginal integration method.

\vskip 4mm

\section{Introduction and motivations}
The multivariate regression function estimation is an important
problem which has been extensively treated for discrete time
processes. It is well-known from \cite{Stone} that the additive
regression models bring out a solution to the problem of the curse
of dimensionality in nonparametric multivariate regression
estimation, which is characterized by a loss in the rate of
convergence of the regression function estimator when the dimension
of the covariates increases. Additive models allow to
reach even univariate rate when these models fit well.\\
 For continuous time processes, \cite{Bosq1993} obtained the optimal rate for the estimator of
  multivariate regression,  which is the same as in the i.i.d.
case. He even proved that, for processes with irregular paths, it is
possible to reach the {\it parametric rate}. This one, called the
superoptimal rate,  does not depend on the dimension of the
variables, but the needed conditions on the processes  are very
strong. That is the reason why it is relevant to study additive
models to bring out a solution to the problem of the curse of
dimensionality.\\

Let ${\bf Z}_t=({\bf X}_t, Y_t), (t\in \mathbb{R})$ be a
$\mathbb{R}^d\times \mathbb{R}$-valued measurable stochastic process
defined on a probability space $(\Omega, \mathcal{A}, P)$. Denote by
$\psi$ a given real measurable function. We consider the additive
regression function associated to $m_\psi(Y)$ defined by,
\begin{eqnarray}\label{add_component}
m_{\psi} ({\bf x})&= & E\left(\psi (Y) \mid \bf{X}=\bf{x}\right),~
\forall ~{\bf x}=(x_1,...,x_d) \in \mathbb{R}^d, \label{fdereg}
\\
&=& \mu + \sum_{l=1}^d m_l(x_l) := m_{\psi,add}(\bf{x}) .
\label{additive}
\end{eqnarray}
\vskip5pt\noindent Let $K_1$, $K_2$, $K_3$ and $K$, be kernels
respectively defined on $\mathbb{R}$, $\mathbb{R}^{d-1}$,
$\mathbb{R}^d$ and $\mathbb{R}^d$. We denote by $\hat f_T$ the
estimate of $f$, the density function of the covariable ${\bf X }$,
(see \cite{Banon1978}), that is,
\begin{eqnarray*}
\hat f_T({\bf x})= \frac{1}{Th_T^d} \int_0^T K\Big(\frac{{\bf
x}-{\bf X}_s}{h_T}\Big) ds,
\end{eqnarray*}
where $(h_T)$ is a positive real function. In estimating the
regression function defined in (\ref{fdereg}), we use the following
two estimators (see for exemple \cite{Bosq1996} and
 \cite{Jones1994})
\begin{equation}
\widetilde{m}_{\psi,T}({\bf x}) = \int_0^T W_{T,t}({\bf
x})\psi(Y_t)dt ~~ \mbox{with}~~W_{T,t}({\bf x})= \frac{
K_{3}\big({\frac{{\bf x} -{\bf
X}_{t}}{h_{1,T}}}\big)}{Th_{1,T}^d\hat f_T ({\bf
X}_t)},\label{estcarbo1}
\end{equation}
and
\begin{equation}
\widetilde{m}_{\psi,T,l}({\bf x}) := \int_0^T W_{T,t}^l({\bf x})
\psi(Y_t) dt~~ \mbox{with}~~W_{T,t}^l({\bf
x})=\frac{K_1\big(\frac{x_l-X_{t,l}}{h_{1,T}}\big)
K_2\big(\frac{{\bf x}_{-l}-{\bf
X}_{t,-l}}{h_{2,T}}\big)}{Th_{1,T}h_{2,T}^{d-1} \hat f_T ({\bf
X}_{t})},\label{estcarbo}
\end{equation}
where  $(h_{j,T}), j=1,2$ are  positive real functions. Let
$q_1,...,q_d$ be $d$ density functions defined in $\mathbb{R}$.
Setting $q({\bf x}) = \prod_{l=1}^d q_l(x_l) ~~\mbox{and}~~
q_{-l}({\bf x}_{-l}) = \prod_{j \neq l} q_j(x_j)$. To estimate the
additive components of the regression function, we use the marginal
integration method (see \cite{Linton1995} and \cite{newey1994}). We
obtain then
\begin{eqnarray}
\eta_{l}(x_{l}) = \int_{\mathbb{R}^{d-1}} m_{\psi}({\bf x})
q_{-l}({\bf x}_{-l}) d{\bf x}_{-l} - \int_{\mathbb{R}^d}
m_{\psi}({\bf x}) q({\bf x}) d{\bf x},\quad l=1,...,d,
\label{additive_component}
\end{eqnarray}
in such a way that the  following two equalities hold,
\begin{eqnarray}
& &\eta_l(x_l) = m_l(x_l) - \int_{\mathbb{R}} m_l(z)
q_l(z)dz, \quad l=1,...,d, \label{relation_additive_component}\\
&&m_{\psi}({\bf x})= \sum_{l=1}^d \eta_l(x_l) +
\int_{\mathbb{R}^d} m_{\psi}({\bf z})q({\bf z}) d{\bf z}
\label{additive_component_marginale}.
\end{eqnarray}
In view of (\ref{relation_additive_component}) and
(\ref{additive_component_marginale}), we note that$\eta_l$ and $m_l$
are equal up to an additional constant. Therefore, $\eta_l$ is also
an additive component, fulfilling a different  identifiability
condition. From (\ref{estcarbo}) and (\ref{additive_component}), a
natural estimate of this $l$-th component is given by
\begin{eqnarray}
\widehat \eta_{l}(x_{l}) = \int_{\mathbb{R}^{d-1}}\!
\widetilde{m}_{\psi,T,l}({\bf x}) q_{-l}({\bf x}_{-l}) d{\bf
x}_{-l} - \int_{\mathbb{R}^d}\! \widetilde{m}_{\psi,T,l}({\bf x})
q({\bf x}) d{\bf x},\ l=1,...,d,\label{etadef}
\end{eqnarray}
from which we deduce the estimate $\widehat m_{\psi,T,add}$ of the
additive regression function,
\begin{eqnarray}
\widehat m_{\psi,T,add}({\bf x}) & = & \sum_{l=1}^d \widehat
\eta_l(x_l) + \int_{\mathbb{R}^d} \widetilde{m}_{\psi,T}({\bf x})
q({\bf x})d{\bf x}. \label{estim_add}
\end{eqnarray}
Before stating our results, we introduce some additional notations
and our  assumptions. Let $\mathcal{C}_1, ...,~ \mathcal{C}_d,$ be
$d$ compact intervals of $\mathbb{R}$ and set
$\mathcal{C}=\mathcal{C}_1\times...\times \mathcal{C}_d$. For every
subset $\mathcal{E}$ of $\mathbb{R}^q$, $q\geq1$, and any
$\delta>0$, introduce the $\delta$-neighborhood $\mathcal{E}^\delta$
of $\mathcal{E}$, namely, $\mathcal{E}^\delta = \{{\bf x} :
\inf_{{\bf y}\in \mathcal{E}}\|{\bf x}-{\bf y}\|_{\mathbb{R}^q}<
\delta\},$ with $\|\cdot\|_{\mathbb{R}^q}$ standing for the
euclidian norm on $\mathbb{R}^q$.\\
\\
(C.1)~~ $\mbox{There exists a positive constant }M \mbox{ such that
$|\psi(y)| \leq M < \infty$}.$\\
\\
(C.2)~~The function $m_\psi \mbox{ is  $k$-times continuously
differentiable, $ k \geq 1$, and }$\\
$$\sup_{{\bf x}} \Big|\frac{\partial^k m_\psi}{\partial
x_\ell^k}({\bf x}) \Big|<\infty;~\ell=1,...,d.$$ Denote by $f_\ell$,
$\ell=1,...,d$ the density functions of $X_\ell$, $\ell=1,...,d$.
The functions $f$ and $f_\ell, \ell=1,...,d,$ are supposed to be
continuous, bounded and \\
\\
$(F.1) \forall{\bf x}\in\mathcal{C}^\delta, ~ f({\bf
x})>0~\mbox{and}~f_{\ell}(x_\ell)>0,\
\ell=1, ..., d.$  \\
\\
$(F.2)~ f \mbox{ is $k'$-times continuously differentiable on
} \mathcal{C}^{\delta}, k'> k d.$  \\
\\
$(F.3)~ \mbox{For all} ~~ 0 < \lambda \leq 1, \Big|\frac{\partial
f^{(k)}}{\partial x_1^{j_1}...\partial_d^{j_d}}({\bf
x}')-\frac{\partial f^{(k)}}{\partial
x_1^{j_1}...\partial_d^{j_d}}({\bf x})\Big| \leq L \|{\bf
x}'-{\bf x}\|^\lambda~\mbox{with}~j_1+...+j_d=k'.$\\
Where $\|.\|$ is a norm on $\mathbb{R}^d$ and $L$ is a positive
constant.
\\

\noindent The kernels $K_1$, $K_2$, $K_3$ and $K$  are assumed  to fulfill the following conditions\\
\\
(K.1)~~  $K_1$, $K_2$, $K_3$ and $K$  are continuous respectively on the compact supports $S_1\subset\mathcal{C}_1$, $S_2\subset\mathcal{C}_2\times...\times\mathcal{C}_d$, $S_3\subset\mathcal{C}$ and $S$,\\
\\
(K.2)~~  $\int K=1$ and $\int K_j=1, ~~j=1,2,3,$\\
\\
(K.3)~~  $\mbox{$K_1$, $K_2$ and $K_3$ are of order $k$,}$\\
\\
(K.4)~~  $\mbox{$K$ is of order $k'$.}$\\
\\
(K.5)~~  $K_1$ is a Lipschitz function. \
\\

\noindent The  density functions  $q_\ell$, $\ell=1, ..., d$,
satisfy
the following assumption\\
\\
$(Q.1)$~~ For any $1\leq l\leq d$,~$q_\ell$ has  k  continuous and
bounded derivatives, with a compact support included in
$\mathcal{C}_\ell.$\\
\\
There exists a set $\Gamma \in \mathcal{B}_{\mathbb{R}^2}$
containing
$D=\{(s,t)\in \mathbb{R}^2: s=t\}$ such that\\
\\
$(D.1)~~  f_{({\bf X}_s,Y_s),({\bf X}_t,Y_t)} -
f_{({\bf X}_s,Y_s)}\bigotimes f_{({\bf X}_t,Y_t)}~~\mbox{exists everywhere for} ~~ (s,t) \in \Gamma^C,$ \\
\\
$(D.2)~~  A_f(\Gamma):= \sup_{(s,t) \in \Gamma^C}\sup_{{\bf
x,y}\in\mathcal{C^{\delta}}\times\mathcal{C^{\delta}}}\int_{u,v\in\mathbb{R}^2}
| f_{({\bf X}_s,Y_s),({\bf X}_t,Y_t)}({\bf x},u,{\bf y},v) \\-
f_{({\bf X}_s,Y_s)}({\bf x},u)
f_{({\bf X}_t,Y_t)}({\bf y},t)|dudv <\infty,$\\
\\
$(D.3)~~\mbox{there exists} ~\ell_\Gamma<\infty~\mbox{and} ~T_0~
\mbox{such that},~ \forall T>T_0,~ \frac{1}{T}
\int_{[0,T]^2\cap \Gamma} ds dt\leq\ell_\Gamma.$ \\
\\
We  work under the following conditions upon the smoothing
parameters $h_T$ and $h_{j,T},~j=1,2$,\\
\\
$(H.1)~~ h_T =c' \Big(\frac{\log T}{T}\Big)^{1/(2k'+d)}, \mbox{ for a fixed } 0<c'<\infty$,\\
\\
$(H.2)~~ h_{1,T}=c_1T^{-1/(2k+1)} ~\mbox{and}~
h_{2,T}=c_2T^{-1/(2k+1)}, \mbox{ for fixed } 0<c_1,c_2<\infty$,\\
\\
$(H.2)'~~h_{1,T}=c_1\big(\frac{\log(T)}{T}\big)^{1/(2k+1)}
~\mbox{and}~
h_{2,T}=c_2\big(\frac{\log(T)}{T}\big)^{1/(2k+1)}, \mbox{ for fixed } 0<c_1,c_2<\infty.$\\
\\
Throughout this work, we use the $\alpha$-mixing dependance
structure where the associated coefficient is defined, for every
$\sigma$-fields $\mathcal{A}$ and $\mathcal{B}$ by
$$\alpha\big(\mathcal{A},\mathcal{B}\big)=\sup_{(A,B)\in(\mathcal{A},\mathcal{B})}|P(A\cap
B)-P(A)P(B)|.$$ For all Borelian set $I$ in $\mathbb{R}^+$ the
$\sigma$-algebra defined by $\big(Z_t,t\in I\big)$ is denoted by
$\sigma\big(Z_t,t\in I\big)$. Writing
$\alpha(u)=\sup_{t\in\mathbb{R}_+}\alpha\big(\sigma\big(Z_v,v\leq
t\big),\sigma\big(Z_v,v\geq t+u\big)\big)$, we use the
condition \\
\\
$(A.1)~~\alpha(t)=\mathcal{O}\big(t^{-b}\big) ~~ \mbox{with}~~b >
2d+10+ \frac{6+4d}{k}$.
\begin{theorem}\label{th1} Under the
conditions $(A.1)$, $(C.1)-(C.2)$, $(F.1)-(F.3)$, $(K.1)-(K.4)$,
$(Q.1)$, $(D.1)-(D.3)$ and $(H.1)-(H.2)$, we have, for  all ${\bf
x}\in\mathcal{C}^\delta$ $$ E(\widehat m_{\psi,T,add}({\bf x}) -
m_\psi({\bf x}))^2 =\mathcal{O}( T^{-2k/2k+1}).$$
\end{theorem}
\begin{theorem}\label{th2}
Under the conditions $(A.1)$, $(C.1)-(C.2)$, $(F.1)-(F.3)$,
$(K.1)-(K.5)$, $(Q.1)$, $(D.1)-(D.3)$ and $(H.1)-(H.2)'$, we have
$$\sup_{{\bf x}\in \mathcal{C}}
|\widehat m_{\psi,T,add}({\bf x}) - m_\psi({\bf x})| =
\mathcal{O}\Bigg(\Big(\frac{\log
T}{T}\Big)^{k/2k+1}\Bigg)~~~\mbox{a.s.}$$
\end{theorem}
\section{Proofs}
The proofs of our theorems are split into two steps. First, we
consider the case where the density is assumed to be known.
Subsequently, we treat the general case  when $f$ is unknown. Denote
by $\hat{\hat{\eta}}$, $\widetilde{\widetilde{m}}_{\psi,T}({\bf x})$
and $\widetilde{\widetilde{m}}_{\psi,T,l}({\bf x})$ the versions of
$\hat{\eta}$, $\widetilde{m}_{\psi,T}({\bf x})$ and
$\widetilde{m}_{\psi,T,l}({\bf x})$ associated to a known (formally,
we replace $\hat{f_T}$  by $f$ in the expressions (\ref{estcarbo1}),
(\ref{estcarbo}) and $\widetilde{m}_{\psi,T,l}({\bf x})$ by
$\widetilde{\widetilde{m}}_{\psi,T,l}({\bf x})$ in (\ref{etadef}).\\
\\
Introduce now the following quantities (see, for the discrete case
\cite{Sarad2000}), we  establish the proof for the first component,
\begin{eqnarray}
&&\widehat{\widehat{m}}_{\psi,T,add}({\bf x}) =  \sum_{l=1}^d
\widehat{\widehat{\eta_l}}(x_l) + \int_{\mathbb{R}^d}
\widetilde{\widetilde{m}}_{\psi,T}({\bf x}) q({\bf x})d{\bf
x}.\label{estim_add} \\
&&\widetilde Y_{\psi,T,t}= \psi (Y_t)
\int_{\mathbb{R}^{d-1}} \frac{1}{h^{d-1}_{2,T}}
K_{2}\Big({\frac{{\bf x}_{-1} - {\bf X}_{t,-1}}{h_{2,T}}}\Big)
\frac{q_{-1}({\bf x}_{-1})}{f(X_{t,-1}|X_{t,1})} d{\bf x}_{-1}, \label{def_Psi_tilde}\\
&&\mathcal{G}({\bf u}_{-1})=\int_{\mathbb{R}^{d-1}}
\frac{1}{h^{d-1}_{2,T}} K_{2}\Big({\frac{{\bf x}_{-1} - {\bf
u}_{-1}}{h_{2,T}}}\Big)
q_{-1}({\bf x}_{-1}) d{\bf x}_{-1}, \label{def_G_cal}\\
&&\widehat \alpha_1(x_1) = \frac{1}{Th_{1,T}} \int_0^T
\frac{\widetilde Y_{\psi,T,t}}{f_1(X_{t,1})} K_{1}
\Big(\frac{x_{1}- X_{t,1}}{h_{1,T}}\Big)dt, ~\mbox{for}~x_1 \in \mathcal{C}_1, \label{def_alpha1}\\
&&\widetilde m_{T}(x_1)= E\big( \left. \widetilde
Y_{\psi,T,t}\right|
X_{t,1} = x_1\big), \label{def_mtilde}\\
&&C_T =  \mu + \int_{\mathbb{R}^{d-1}} \sum_{j=2}^d m_j(u_j)
\mathcal{G}({\bf u}_{-1}) d{\bf u}_{-1},\label{def_Cn}\\
&&\widehat C_T=\int_{\mathbb{R}^d} \widetilde{\widetilde
m}_{\psi,T,1}({\bf x}) q({\bf x})d{\bf x}, \label{def_Cn_hat}\\
&&C=\int_{\mathbb{R}}m_1(x_1)q_1(x_1)dx_1.
\end{eqnarray}
The following Lemma  is of particular interest to establish the
result of theorem (\ref{th1}). Note that (\ref{equation3}) is
``only" be instrumental in the proof of (\ref{equation1}).
\begin{lemma}\label{lem1}
Under the assumptions $(C.1)-(C.2)$, $(F.1)-(F.2)$, $(K.1)$,
$(Q.1)$ and $ (H .2)$ , we have
\begin{eqnarray}
&&E(\hat C_T -C_T+C)^2
=\mathcal{O}\Big(T^{-2k/(2k+1)}\Big),\label{equation2}\\
&&{\rm Var}(\hat
\alpha_1(x_1))=\mathcal{O}\Big(T^{-2k/(2k+1)}\Big),
\label{equation3}\\
& &E\big(\widehat {\widehat \eta}_1(x_1) - \eta_1(x_1)\big)^2 =
\mathcal{O}\Big(T^{-2k/(2k+1)}\Big).
 \label{equation1}
\end{eqnarray}
\end{lemma}
{\it Proof:}~~ According to Fubini's Theorem and under the additive
model assumption, we have
\begin{eqnarray*}
\mathbb{E}(\hat C_T -C_T) & = & \mathbb{E}\Big\{\int_{\mathbb{R}^d}
\widetilde {\widetilde m}_{\psi,T,1}({\bf x}) q({\bf x})d{\bf x} -
\mu - \int_{\mathbb{R}^{d-1}} \sum_{j=2}^d m_j(u_j)\mathcal{G}({\bf
u}_{-1})d{\bf u}_{-1}\Big\}\\
& =&\int_{\mathbb{R}^d}E( \widetilde {\widetilde m}_{\psi,T,1}({\bf
x})) q({\bf x})d{\bf x}  - \mu - \int_{\mathbb{R}^{d-1}}
\sum_{j=2}^d m_j(u_j)\mathcal{G}({\bf u}_{-1})d{\bf u}_{-1}\\
& = & \int_{\mathbb{R}^d}\frac{1}{h_{1,t}}m_{\psi}({\bf
u})\mathcal{G}({\bf u}_{-1})\int_{\mathbb{R}}K_1\Big(\frac{x_1-
u_1}{h_{1,T}}\Big)q_1(x_1)dx_1 d{\bf u} \\ &&
~~~~~~~~~~~~~~~~~~~~~~~~~~~~~~~~~~~~~~~~~~- \mu -
\int_{\mathbb{R}^{d-1}} \sum_{j=2}^d m_j(u_j)
\mathcal{G}({\bf u}_{-1}) d{\bf u}_{-1}\\
& = & \sum_{j=1}^d \int_{\mathbb{R}^d}
\frac{1}{h_{1,T}}m_{j}(u_j)\mathcal{G}({\bf
u}_{-1})\int_{\mathbb{R}}K_1\Big(\frac{x_1-
u_1}{h_{1,T}}\Big)q_1(x_1)dx_1 d{\bf u} \\&&
~~~~~~~~~~~~~~~~~~~~~~~~~~~~~~~~~~~~~~~~~~~~~~~~-\int_{\mathbb{R}^{d-1}}
\sum_{j=2}^d m_j(u_j) \mathcal{G}({\bf
u}_{-1}) d{\bf u}_{-1}\\
&  = &\int_{\mathbb{R}}\int_{\mathbb{R}}
\frac{1}{h_{1,T}}m_1(u_1)K_1\Big(\frac{x_1-
u_1}{h_{1,T}}\Big)q_1(x_1) dx_1du_1.
\end{eqnarray*}
Setting $v_1 h_{1,T}=x_1- u_1$ and using a Taylor expansion, we get,
by $(C.2)$ and $(K.1)-(K.3)$,
\begin{eqnarray}\label{ii}
\mathbb{E}(\hat C_T -C_T) - C \nonumber
& = & \int_{\mathbb{R}}\int_{\mathbb{R}}q_1(x_1)m_1(x_1-h_{1,T}v_1)K_1(v_1)dv_1dx_1 - C\nonumber\\
& = & \int_{\mathbb{R}}\int_{\mathbb{R}}q_1(x_1)
[m_1(x_1-h_{1,T}v_1)-m_1(x_1)]K_1(v_1)dv_1dx_1\nonumber\\
& = &
\int_{\mathbb{R}}\int_{\mathbb{R}}q_1(x_1)\Big[\frac{(-h_{1,T})^k
v_1^k}{k!} m_1^{(k)}(x_1)\Big]K_1(v_1)dv_1dx_1 \\
&& + o(h_{1,T}^k).\nonumber
\end{eqnarray}
Under $(H.2)$, it follows that,
\begin{eqnarray}
\Big[\mathbb{E}\Big(\hat C_T -C_T - C\Big)\Big]^2=
\mathcal{O}\big(T^{-2k/(2k+1)}\big). \label{resulA}
\end{eqnarray}
The Fubini's theorem gives us
\begin{eqnarray}
{\rm Var} (\hat C_T) &=& \frac{1}{(Th_{1,T})^2}{\rm Var} \Bigg(
\int_{0}^T\frac{\psi( Y_{t})}{f({\bf X}_t)}{\mathcal G}({\bf
X}_{t, -1})\int_{\mathbb{R}}K_1\Big(\frac{x_1-X_{t,1}}
{h_{1,T}}\Big)q_1(x_1)dx_1 dt\Bigg)\nonumber\\
&=& \frac{1}{(Th_{1,T})^2} \int_{{t,s}\in[0;T]^2}{\rm
Cov}\Bigg(\frac{\psi (Y_{t})}{f({\bf X}_t)}{\mathcal G}({\bf
X}_{t,-1})\int_{\mathbb{R}}K_1\Big(\frac{x_1-X_{t,1}}
{h_{1,T}}\Big)\nonumber \\
&&q_1(x_1)dx_1;\frac{\psi( Y_s)}{f({\bf X}_s)}{\mathcal G}({\bf
X}_{s, -1})\int_{\mathbb{R}}K_1\Big(\frac{y_1-X_{s,1}}
{h_{1,T}}\Big)q_1(y_1)dy_1\Bigg)dsdt.\end{eqnarray} Under $(C.1)$,
$(F.1)$, $(K.1)-(K.2)$ and $(Q.1)$, there exists a finite constant
$M_3$ such that, for T large enough,
$$\inf\Bigg\{a~:~P\Big(\frac{\psi (Y_{t})}{f({\bf X}_t)}{\mathcal G}({\bf
X}_{t, -1})\int_{\mathbb{R}}K_1\Big(\frac{x_1-X_{t,1}}
{h_{1,T}}\Big)q_1(x_1)dx_1>a\Bigg)=0\Big\}\leq h_{1,T}M_3.$$ Thus,
using the Billingsley's inequality and the condition $(A.1)$,
\begin{eqnarray} {\rm Var} (\hat C_T)&\leq &
\frac{8M_3^2}{T^2}\int_{s\in[0;T]}\int_{t\in[0,T-s]}\alpha_{t}dtds
=\mathcal{O}\Big(\frac{1}{T}\Big).\label{resulB}
\end{eqnarray}
Finally, by combining the statements  $(\ref{resulA})$ and
$(\ref{resulB})$, we obtain (\ref{equation2}).
\vskip5pt \noindent
\emph{Proof of (\ref{equation3}).} Recalling (\ref{def_alpha1}),
we have
\begin{eqnarray}
{\rm Var}( \hat \alpha_{1}(x_{1}))& = &
\frac{1}{T^2}\bigg[\int_{[0,T]^2\cap\Gamma} {\rm Cov}\Big(\frac{\widetilde Y_{\psi,T,t}}{f_{1}(X_{t,1})h_{1,T}}K_1\Big(\frac{x_1-X_{t,1}}{h_{1,T}}\Big),\nonumber\\
&& \frac{\widetilde
Y_{\psi,T,s}}{f_{1}(X_{s,1})h_{1,T}}K_1\Big(\frac{x_1-X_{s,1}}{h_{1,T}}\Big)\Big)dsdt
\nonumber \\
&& +\int_{[0,T]^2\cap\Gamma^c}{\rm Cov}\Big(\frac{\widetilde Y_{\psi,T,t}}{f_{1}(X_{i,1})h_{1,T}}K_1\Big(\frac{x_1-X_{t,1}}{h_{1,T}}\Big),\nonumber\\
&&\frac{\widetilde
Y_{\psi,T,s}}{f_{1}(X_{s,1})h_{1,T}}K_1\Big(\frac{x_1-X_{s,1}}{h_{1,T}}\Big)\Big)dsdt
  \bigg]\nonumber~\\&:=&A+B.\label{chepa}
\end{eqnarray}
For the first term, noting that, under $(C.1)$, $(F.1)$,
$(K.1)-(K.2)$ and $(Q.1)$, there exists a finite constant $M_4$
such that, for $T$ large enough,
$$M_4\geq\inf\Bigg\{a~:~P\Bigg(\frac{\widetilde Y_{\psi,T,0}^2}{f_1( X_{0,1})^2}\Big|K_1\left(\frac{x_1-X_{0,1}}{h_{1,T}}\right)\Big|>a\Bigg)=0\Bigg\}.$$
Thus, we have
\begin{eqnarray}A &\leq& \frac{1}{T^2h_{1,T}^{2}}\int_{[0,T]^2\cap \Gamma
}\mathbb{E} \left(\frac{\widetilde
Y_{\psi,T,0}}{f_{1}(X_{0,1})}K_1\left(\frac{x_1-X_{0,1}}{h_{1,T}}\right)\right)^2dsdt\nonumber
\\&\leq& \frac{M_4}{T^2h_{1,T}^{2}}\int_{[0,T]^2\cap
\Gamma }\int_{\mathbb{R}}
\Big|K_1\left(\frac{x_1-u}{h_{1,T}}\right)\Big||f_1(u)|dudsdt\nonumber
\\&\leq& \frac{M_4\|f\|_{\infty}l_{\Gamma}}{Th_{1,T}}\int_{\mathbb{R}}
\big|K_1\left(v\right)\big|dv
=\mathcal{O}\left(\frac{1}{Th_{1,T}}\right)\label{A}.\end{eqnarray}
To treat the second term, we introduce the set
$S_{a(T)}=\{(s,t)\in\mathbb{R}^2;|t-s|\leq a(T)\}$, where $a(T) =
h_T^{-1}$, we have
\begin{eqnarray} B &=&\frac{1}{T^2}\int_{[0,T]^2\cap\Gamma^c\cap S_{a(T)}}{\rm Cov}\bigg(\frac{\widetilde Y_{\psi,T,t}}{f_{1}(X_{t,1})h_{1,T}}K_1\Big(\frac{x_1-X_t}{h_{1,T}}\Big),\nonumber\\
&&\frac{\widetilde
Y_{\psi,T,s}}{f_{1}(X_{s,1})h_{1,T}}K_1\Big(\frac{x_1-X_s}{h_{1,T}}\Big)\bigg)dsdt
  \nonumber \\ &+&\frac{1}{T^2}\int_{[0,T]^2\cap\Gamma^c\cap S_{a(T)}^c}{\rm Cov}\bigg(\frac{\widetilde Y_{\psi,T,t}}{f_{1}(X_{t,1})h_{1,T}}K_1\Big(\frac{x_1-X_t}{h_{1,T}}\Big),\nonumber\\
&&\frac{\widetilde
Y_{T,s}}{f_{1}(X_{s,1})h_{1,T}}K_1\Big(\frac{x_1-X_s}{h_{1,T}}\Big)\bigg)dsdt
  \nonumber \\
  &:=&E+F \label{B}.\end{eqnarray}
Under the conditions $(C.1)$, $(F.1)$, $(K.1)-(K.2)$ and $(Q.1)$,
there exists a constant $M_5$ such that, for $T$ large enough,
$$\sup_{z_1\in\mathbb{R}}\int_{(y,{\bf z_{-1}})\in\mathbb{R}\times \mathbb{R}^{d-1}}\Bigg|\frac{\psi (y)\mathcal{G}({\bf z}_{-1})}{f({\bf z})} \mathds{1}_{S_1}\Big(\frac{x_1-z_1}{h_{1,T}}\Big)\Bigg|dyd{\bf z_{-1}}\leq M_5.$$
Consider now the term $E$, we have
\begin{eqnarray}E&=& \frac{1}{T^2h_{1,T}^{2}}\int_{(s,t)\in[0,T]^2\cap\Gamma^c\cap S_{a(T)}}\int_{(u,{\bf v})\in\mathbb{R}\times\mathbb{R}^d}\int_{(y,{\bf z})\in\mathbb{R}\times\mathbb{R}^d}\frac{\psi (y)\mathcal{G}({\bf z}_{-1})}{f({\bf z})}\nonumber
\\
&&K_1\Big(\frac{x_1-z_1}{h_{1,T}}\Big)\frac{\psi (u)\mathcal{G}({\bf v}_{-1})}{f({\bf v})}K_1\Big(\frac{x_1-v_1}{h_{1,T}}\Big)\Big(f_{Y_t,{\bf X}_{t},Y_s,{\bf X}_{s}}(u,{\bf v},y,{\bf z})\nonumber \\ &-&f_{Y_t,{bf X}_{t}}(u,{\bf v})f_{Y_s,{bf X}_{s}}(y,{\bf z})\Big)dyd{\bf z}dud{\bf v}dsdt\nonumber \\
&\leq&\frac{2a(T)M_5^2\|K_1\|_{\mathbb{L}_1}^2A_f(\Gamma)}{T}.\label{m1}
\end{eqnarray}
Noting that, under the conditions $(C.1)$, $(F.1)$, $(K.1)-(K.2)$
and $(Q.1)$, there exists a finite constant $M_6$ such that, for
$T$ large enough,
$$M_6\geq\inf\Bigg\{a~:~P\Bigg(\frac{\widetilde Y_{\psi,T,0}}{f_1(X_{0,1})}\Big|K_1\left(\frac{x_1-X_{0,1}}{h_{1,T}}\right)\Big|>a\Bigg)=0\Bigg\}.$$
Using the Billingsley's inequality, it follows that
 \begin{eqnarray}F&\leq& \frac{2}{T^2h_{1,T}^{2}} \int_{[0,T]^2\cap \Gamma^c\cap S_{a(T)}^c\cap\{u>v\}}4M_6^2 \alpha(u-v)dudv\nonumber \\
&\leq&
\frac{8M_6^2}{Th_{1,T}^{2}}\int_{\{t>a(T)\}}\alpha(t)dt\nonumber\\
&\leq& \frac{8M_6^2}{Th_{1,T}^{2}}La(T)^{-1}.\label{m2}
\end{eqnarray}
Finally, combining the hypothesis $(H.2)$ and  the statements
(\ref{chepa})and (\ref{m2}), we obtain
(\ref{equation3}).\\
\noindent \emph{Proof of (\ref{equation1})}. We have
\begin{eqnarray}
&& \hspace{-1.5cm}\mathbb{E}(\hat \alpha_1(x_1) ) - \widetilde m_{T}(x_1)\nonumber\\
&=& \int_{\mathbb{R}}\frac{1}{h_{1,T}}
\widetilde m_{T}(u_1) K_1\Big( \frac{x_1-u_1}{h_{1,T}} \Big) du_1 -\widetilde m_{T}(x_1)\nonumber\\
&=& \int_{\mathbb{R}} \big[\widetilde
m_{T}(x_1-v_1h_{1,T})-\widetilde
m_{T}(x_1)\big] K_1(v_1)dv_1\nonumber\\
&=&\int_{\mathbb{R}} \int_{\mathbb{R}^{d-1}}
[m_\psi(x_1-v_1h_{1,T},{\bf u}_{-1}) - m_\psi(x_1,{\bf
u}_{-1})]\mathcal{G}({\bf u}_{-1}) d{\bf u}_{-1} K_1(v_1)dv_1\nonumber\\
& = & \int_{\mathbb{R}} \int_{\mathbb{R}^{d-1}}
\left[\frac{(-h_{1,T}v_1)^k}{k!} \frac{\partial^k
m_{\psi}}{\partial v_1^k}( v_1,{\bf u}_{-1} )\right]
\mathcal{G}({\bf u}_{-1}) d{\bf u}_{-1}  K_1(v_1)dv_1 +
o(h_{1,T}^k).\nonumber \label{controle1}
\end{eqnarray}
Under the condition $(H.2)$, we obtain
\begin{eqnarray}
\Big[\mathbb{E}\Big(\hat \alpha_1(x_1)- \widetilde
m_{T}(x_1)\Big)\Big]^2=
\mathcal{O}\Big(T^{-2k/(2k+1)}\Big).\label{intermed1}
\end{eqnarray}
Thus, we have
\begin{eqnarray}
\mathbb{E}\Big(\hat \alpha_1(x_1)- \widetilde
m_{T}(x_1)\Big)^2&=&\Big[\mathbb{E}\Big(\hat \alpha_1(x_1)-
\widetilde m_{T}(x_1)\Big)\Big]^2+{\rm Var}(\hat
\alpha_1(x_1)).\label{intermed}
\end{eqnarray}
Consequently, by combining the following inequality
\begin{eqnarray}
E\big(\widehat {\widehat \eta}_1(x_1) - \eta_1(x_1)\big)^2  \leq
2\mathbb{E}\Big(\hat \alpha_1(x_1)- \widetilde m_{T}(x_1)\Big)^2 + 2
E(\hat C_T - C_T -C)^2,  \label{composante}
\end{eqnarray}
and the statements (\ref{intermed1}),  (\ref{intermed}),
(\ref{equation3}) and (\ref{composante}), the proof of
(\ref{equation1}) is readily achieved.
\vskip10pt
\subsection{Proof of Theorem \ref{th1}}
\noindent Using the classical inequality
$(a+b)^2\leq 2(a^2+b^2)$, if follows that, for all ${\bf
x}\in\mathcal{C}$,
\begin{eqnarray}
E(\widehat m_{\psi,T,add}({\bf x}) - m_\psi({\bf x}))^2 & \leq &
2E\big( \widehat {\widehat m}_{\psi,T,add}({\bf x})- m_\psi ({\bf
x})\big)^2+ 2 E\big(\widehat m_{\psi,T,add}({\bf x})- \widehat
{\widehat m}_{\psi,T,add}({\bf x}) \big)^2\nonumber \\
&:=& I_1({\bf x})+I_2({\bf x}).
\end{eqnarray}
First, consider the term $I_1$, we have
 \begin{eqnarray}\label{I1majoration} I_1({\bf
x})&=&2E\big( \widehat {\widehat m}_{\psi,T,add}({\bf x})- m_\psi
({\bf x})\big)^2\nonumber\\
&\leq &4d \sum_{\ell=1}^d E\big(\widehat {\widehat
\eta}_\ell(x_\ell) - \eta_\ell(x_\ell)\big)^2 + 4
E\Big[\int_{\mathbb{R}^d}\big(\widetilde {\widetilde
m}_{\psi,T}({\bf x})- m_{\psi,T}({\bf x})\big) q({\bf x})d{\bf
x}\Big]^2.
\end{eqnarray}
Arguing as in proof of Lemma (\ref{lem1}), we obtain
\begin{eqnarray*}
\widehat C_T - C_T -C \!\!\!\! &=&\!\!\!\! \int_{\mathbb{R}^d}\!\!
\widetilde {\widetilde m}_{\psi,T}({\bf x})q({\bf x})d{\bf x}\!\!  -
\!\! \mu\!\!  -\!\! \int_{\mathbb{R}^{d-1}} \!\! \sum_{j=2}^d \!\!
m_j(u_j)
\mathcal{G}({\bf u}_{-1}) d{\bf u}_{-1}-\int_{\mathbb{R}}m_1(x_1)q_1(x_1)dx_1,\\
&=&\int_{\mathbb{R}^d}\widetilde {\widetilde m}_{\psi,T}({\bf
x})q({\bf x})d{\bf x} - \mu -\int_{\mathbb{R}^{d-1}} \sum_{j=1}^d
m_j(u_j)
q({\bf u}) d{\bf u}+\mathcal{O}\Big(h_{1,T}^k\Big),\\
&=&\int_{\mathbb{R}^d}\big(\widetilde {\widetilde m}_{\psi,T}({\bf
x})- m_{\psi,T}({\bf x})\big) q({\bf x})d{\bf
x}+\mathcal{O}\Big(T^{-k/(2k+1)}\Big).
\end{eqnarray*}
It follows that,
\begin{eqnarray}\label{BiasI1}
E\Big[\int_{\mathbb{R}^d}\big(\widetilde {\widetilde
m}_{\psi,T}({\bf x})- m_{\psi,T}({\bf x})\big) q({\bf x})d{\bf
x}\Big]^2 &\leq & 2E\Big(\widehat C_T - C_T -C
\Big)^2+\mathcal{O}\Big(T^{-2k/(2k+1)}\Big).
\end{eqnarray}
By combining  (\ref{I1majoration}), (\ref{BiasI1}),
(\ref{equation2}), and (\ref{equation1}), we conclude that, for all
${\bf x}\in \mathcal{C}$
\begin{eqnarray}
I_1({\bf x})=\mathcal{O}\Big(T^{-2k/(2k+1)}\Big). \label{I1_proof}
\end{eqnarray}
\vskip15pt
\noindent  Turning our attention to $I_2({\bf x})$, it
holds that,
\begin{eqnarray*}
& &\hspace{-2cm} E\big(\widehat m_{\psi,T,add}({\bf x})- \widehat
{\widehat
m}_{\psi,T,add}({\bf x}) \big)^2\\
& = & E\left[\sum_{\ell=1}^d (\widehat {\widehat
\eta}_\ell(x_\ell) - \widehat \eta_\ell (x_\ell))+
\int_{\mathbb{R}^d} \widetilde {\widetilde m}_{\psi,T} ({\bf
x})q({\bf x}) d{\bf x} - \int_{\mathbb{R}^d}\widetilde
m_{\psi,T}({\bf x}) q({\bf x}) d{\bf
x}\right]^2\\
& \leq & 4d \sum_{\ell=1}^d E\Big(\int_{\mathbb{R}^{d-1}}
(\widetilde
 m_{\psi,T,\ell} ({\bf x}) -\widetilde {\widetilde
 m}_{\psi,T,\ell} ({\bf x})) q({\bf x}_{-\ell})d{\bf x}_{-\ell}\Big)^2\nonumber \\ & &+ 4d
 \sum_{\ell=1}^d E\Big(\int_{\mathbb{R}^{d}} (\widetilde
 m_{\psi,T,\ell} ({\bf x}) -\widetilde {\widetilde
 m}_{\psi,T,\ell} ({\bf x})) q({\bf x}) d{\bf x}\Big)^2\\
 &  &+ 2E\left(\int_{\mathbb{R}^d} \widetilde {\widetilde
m}_{\psi,T} ({\bf x})q({\bf x}) d{\bf x} -
\int_{\mathbb{R}^d}\widetilde m_{\psi,T}({\bf x}) q({\bf x}) d{\bf
x}\right)^2\\
& \leq & 4d \sum_{\ell=1}^d E\int_{\mathbb{R}^{d-1}} (\widetilde
 m_{\psi,T,\ell} ({\bf x}) -\widetilde {\widetilde
 m}_{\psi,T,\ell} ({\bf x}))^2 q^2({\bf x}_{-\ell}) d{\bf x}_{-\ell} \nonumber \\ & &+ 4d
 \sum_{\ell=1}^d E\int_{\mathbb{R}^{d}} (\widetilde
 m_{\psi,T,\ell} ({\bf x}) -\widetilde {\widetilde
 m}_{\psi,T,\ell} ({\bf x}))^2 q^2({\bf x}) d{\bf x}\\
 &  & 2E\int_{\mathbb{R}^d} \big(\widetilde {\widetilde
m}_{\psi,T} ({\bf x})- \widetilde m_{\psi,T}({\bf x})\big)^2
q^2({\bf x})
d{\bf x}\\
& \leq &4d \sum_{\ell=1}^d \int_{\mathbb{R}^{d-1}} E(\widetilde
 m_{\psi,T,\ell} ({\bf x}) -\widetilde {\widetilde
 m}_{\psi,T,\ell} ({\bf x}))^2 q^2({\bf x}_{-\ell}) d{\bf x}_{-\ell} \nonumber \\ &&+ 4d
 \sum_{\ell=1}^d \int_{\mathbb{R}^{d}} E(\widetilde
 m_{\psi,T,\ell} ({\bf x}) -\widetilde {\widetilde
 m}_{\psi,T,\ell} ({\bf x}))^2 q^2({\bf x}) d{\bf x}\\
 &  & + 2\int_{\mathbb{R}^d} E\big(\widetilde {\widetilde
m}_{\psi,T} ({\bf x})- \widetilde m_{\psi,T}({\bf x})\big)^2
q^2({\bf x}) d{\bf x}
\end{eqnarray*}
Using the decomposition $1/f= 1/\widehat f_T + (\widehat f_T
-f)/(\widehat f_T f)$, it is easily shown that for some positive
constant $M_1<\infty$, we have, under $(Q.2)$, for all ${\bf
x}\in\mathcal{C}$ and $T$ large enough,
\begin{eqnarray*}
&  & E\Big( \widetilde m_{\psi,T,\ell}({\bf x}) - \widetilde
{\widetilde
m}_{\psi,T,\ell}({\bf x})\Big)^2\\
& \leq & M_1
E\left(\frac{1}{h_{1,T}h_{2,T}^{d-1}}\Big|K_1\Big(\frac{x_\ell-X_{t,\ell}}{h_{1,T}}\Big)
K_2\Big( \frac{{\bf x}_{-\ell}- {\bf X}_{t,-\ell}
}{h_{2,T}}\big)\Big|\times \sup_{{\bf x}\in\mathcal{C}}|\widehat
f_T({\bf x})-f({\bf x})|\right)^2 .
\end{eqnarray*}
It's easily seen that under our assumptions, following the
demonstration of Theorem 4.9. in \cite{Bosq1993} p.112 and
replacing $\log_m$ by $1$, we have,
\begin{eqnarray*}
\sup_{{\bf x}\in\mathcal{C}}|\widehat f_T({\bf x})-f({\bf x})|=
\mathcal{O}\Big(\Big(\frac{\log T}
{T}\Big)^{k'/(2k'+d)}\Big)~~\mbox{almost surely },
\end{eqnarray*}
We conclude that, for all ${\bf x} \in \mathcal{C}$,
\begin{eqnarray}
 E\Big( \widehat {\widehat
m}_{\psi,T,add}({\bf x})-\widehat m_{\psi,T,add}({\bf x}) \Big)^2=
\mathcal{O}\Big(\Big(\frac{\log
T}{T}\Big)^{2k'/(2k'+d)}\Big)=\mathcal{O}\Big(T^{-2k/(2k+1)}
\Big).\label{I2_proof}~~~~~~~\sqcap
\end{eqnarray}
\vskip15pt
\subsection{Proof of Theorem \ref{th2}}
 \noindent

 In the next lemma we evaluate the difference between
the estimator of the additive regression function
$\widehat{\widehat{m}}_{\psi,T,add}$,  for continuous time process,
and  the estimator $\widehat{\widehat{m}}_{\psi,n,add}$ where
$n\in\mathbb{N}$.
\begin{lemma}\label{lem3}
For $n\in\mathbb{N}$  large enough, there exists a deterministic
constant C such that for all $\omega $ in $\Omega$ and for all  $T$
in $[n,n+1[$,
\begin{eqnarray*}
|\widehat{\widehat{m}}_{\psi,T,add}({\bf
x})(\omega)-E\widehat{\widehat{m}}_{\psi,T,add}({\bf
x})(\omega)-\widehat{\widehat{m}}_{\psi,n,add}({\bf
x})(\omega)+E\widehat{\widehat{m}}_{\psi,n,add}({\bf
x})(\omega)|<C\bigg(\frac{\log(T)}{T}\bigg)^{\frac{k}{2k+1}}.
\end{eqnarray*}
\end{lemma}
{\it Proof:}~~ It is sufficient to prove that $\forall$
 $\omega\in\Omega,~\forall T\in[n,n+1[,$
$$|\widehat{\widehat{m}}_{\psi,T,add}({\bf
x}(\omega)-\widehat{\widehat{m}}_{\psi,n,add}({\bf
x})(\omega)|<C'\bigg(\frac{\log(T)}{T}\bigg)^{\frac{k}{2k+1}},$$ the
other part being a trivial consequence of this inequality. Moreover,
in view of (\ref{etadef}) and (\ref{estim_add}), we can establish
the following inequalities
\begin{eqnarray}&& \hspace{-1.5cm}\int_{\mathbb{R}^d}\widetilde{\widetilde{m}}_{\psi,T}({\bf x})(\omega)q({\bf x})d{\bf x}-\int_{\mathbb{R}^d}\widetilde{\widetilde{m}}_{\psi,n}({\bf x})(\omega)q({\bf x})d{\bf x}<C_1\bigg(\frac{\log(T)}{T}\bigg)^{\frac{k}{2k+1}},\label{seul}\\
&&\hspace{-1.5cm}\int_{\mathbb{R}^{d-1}}\widetilde{\widetilde{m}}_{\psi,T,l}({\bf
x})(\omega)q_{-l}({\bf x}_{-l})d{\bf
x}_{-l}-\int_{\mathbb{R}^{d-1}}\widetilde{\widetilde{m}}_{\psi,n,l}({\bf
x})(\omega)q_{-l}({\bf x}_{-l})d{\bf
x}_{-l}<C_l'\bigg(\frac{\log(T)}{T}\bigg)^{\frac{k}{2k+1}},\label{trop1}\\
&&\hspace{-1.5cm}\int_{\mathbb{R}^d}\widetilde{\widetilde{m}}_{\psi,T,l}({\bf
x})(\omega)q({\bf x})d{\bf
x}-\int_{\mathbb{R}^d}\widetilde{\widetilde{m}}_{\psi,n,l}({\bf
x})(\omega)q({\bf x})d{\bf
x}<C_l''\bigg(\frac{\log(T)}{T}\bigg)^{\frac{k}{2k+1}},\label{trop2}\end{eqnarray}
with $l=1,...d.$ We  just establish the first inequality, the
techniques being the same for (\ref{trop1}) and (\ref{trop2}). For
fixed  $\omega$ in $\Omega$ and ${\bf x}$ in $\mathbb{R}^d$, we
have, for $n$ large enough,  $$y\notin \mathcal{C}^\delta
\Rightarrow K_3\Big(\frac{{\bf x}-{\bf y}}{h_t}\Big)q({\bf
x})=0,~~~\forall t\geq n.$$ So, by Fubini's Theorem
\begin{eqnarray*}
\int_{\mathbb{R}^d}\widetilde{\widetilde{m}}_{\psi,T}({\bf
x})(\omega)q({\bf x})d{\bf
x}&=&\frac{1}{Th_{1,T}^d}\int_{[0,T]}\int_{\mathbb{R}^d}\frac{\psi(Y_t(\omega))}{f({\bf
X}_t(\omega))}K_3(\frac{{\bf x}-{\bf X}_t(\omega)}{h_{1,T}})q({ \bf
x})d{\bf x}dt\nonumber
\\
&=&\frac{1}{Th_{1,T}^d}\int_{[0,n]}\int_{\mathbb{R}^d}\frac{\psi(Y_t(\omega))}{f({\bf
X}_t(\omega))}K_3(\frac{{\bf x}-{\bf X}_t(\omega)}{h_{1,T}})q({\bf
x})d{\bf x}dt
+ \mathcal{O}(\frac{1}{T})\\
&=&\frac{1}{nh_{1,T}^d}\int_{[0,n]}\int_{\mathbb{R}^d}\frac{\psi(Y_t(\omega))}{f({\bf
X}_t(\omega))}K_3(\frac{{\bf x}-{\bf X}_t(\omega)}{h_{1,T}})q({\bf
x})d{\bf x}dt+\mathcal{O}(\frac{1}{T})
\nonumber\\
&=&\frac{1}{n}\int_{[0,n]}\int_{\mathbb{R}^d}\frac{\psi(Y_t(\omega))}{f({\bf
X}_t(\omega))}K_3({\bf u})q({\bf X}_t(\omega)+{\bf u}h_{1,T})d{\bf
u}dt+\mathcal{O}(\frac{1}{T}).\end{eqnarray*} Denoting
$M_7:=\sup_{{\bf
x}\in\mathcal{C}^\delta,y\in\mathbb{R}}\frac{\psi(y)}{f({\bf
x})}<\infty$, we have
\begin{eqnarray}&& \hspace{-1.5cm}\Big|\int_{\mathbb{R}^d}\widetilde{\widetilde{m}}_{\psi,T}({\bf x})(\omega)q({\bf x})d{\bf x}-\int_{\mathbb{R}^d}\widetilde{\widetilde{m}}_{\psi,n}({\bf x})(\omega)q({\bf x})d{\bf x}\Big|\nonumber\\ &=&\frac{1}{n}\bigg|\int_{t\in[0,n]}\int_{u\in\mathbb{R}^d}\frac{\psi(Y_t(\omega))}{f(X_t(\omega)}K_3({\bf u})(q({\bf X}_t+h_{1,T}{\bf u})-q({\bf X}_t+h_{1,n}{\bf u}))dudt\Big|+\mathcal{O}(\frac{1}{T})\nonumber\\
&\leq&\frac{M_7}{n}\Big|\int_{[0,n]}\int_{\mathbb{R}^d}K_3({\bf u})(q({\bf X}_t(\omega)+h_{1,T}{\bf u})-q({\bf X}_t(\omega)+h_{1,n}{\bf u}))d{\bf u}dt\Big|+\mathcal{O}(\frac{1}{T})\nonumber\\
&\leq&\frac{dM_7}{n}\int_{[0,n]}\int_{
S_3}|K_3({\bf u})\max_{1 \leq l \leq d}\|\frac{\partial q}{\partial u_l}\|_{\infty}\|{\bf u}\|(h_{1,T}-h_{1,n})|d{\bf u}dt+\mathcal{O}(\frac{1}{T})\nonumber\\
&=&\mathcal{O}(h_{1,T}-h_{1,n})+\mathcal{O}(\frac{1}{T})
\end{eqnarray}
Which implies (\ref{seul}) by (K.1). This  achieves the proof of
Lemma \ref{lem3}.\\
\\
Set $\epsilon^2(T)=C\bigg(\frac{\log T}{T}\bigg)^{\frac{2k}{2k+1}},$
where $C$ is a finite constant. There exists a finite number $r(T)
:= \big(3/M'h_{1,T}^2\varepsilon(T)\big)^d$ of balls $B_p$ of center
${\bf x}_p$ and radius $h_{1,T}^2\varepsilon(T)$, such that
$\mathcal{C}\subset \cup_{p=1}^{r(T)} B_p$, where M' is a constant.
For each ${\bf x} \in B_p$ we  denote $t({\bf x})={\bf x}_p$.Write
now,
\begin{eqnarray*}
&  &\hspace{-1.5cm}\sup_{{\bf x} \in \mathcal{C}}|\widehat {\widehat
m}_{\psi,T,add} ({\bf x}) - m_{\psi}({\bf
x})|\\
& \leq & \sup_{{\bf x} \in \mathcal{C}} |\widehat {\widehat
m}_{\psi,T,add} ({\bf x}) - \widehat m_{\psi,T,add} ({\bf x})| +
\sup_{{\bf x} \in \mathcal{C}} |E\widehat {\widehat m}_{\psi,T,add}
({\bf x}) - m_{\psi} ({\bf x})|\\&& +  \sup_{{\bf x} \in
\mathcal{C}}|\widehat {\widehat m}_{\psi,T,add} ({\bf x}) - \widehat
{\widehat m}_{\psi,T,add} (t({\bf x}))|+ \sup_{{\bf x} \in
\mathcal{C}}|E\widehat {\widehat m}_{\psi,T,add} ({\bf x}) -
E\widehat {\widehat m}_{\psi,T,add} (t({\bf x}))|\\&& + \sup_{{\bf
x} \in \mathcal{C}}|\widehat {\widehat m}_{\psi,add} (t({\bf x})) -
E\widehat {\widehat m}_{\psi,add} (t({\bf x}))|.
\end{eqnarray*} Thus, to prove the Theorem \ref{th2}, it
suffices to establish the following Lemma.
\begin{lemma} Under the same hypothesis as Theorem \ref{th2}, we have
\begin{eqnarray}
&  &\sup_{{\bf x} \in \mathcal{C}}|E\widehat {\widehat
m}_{\psi,T,add} ({\bf x}) - m_{\psi}
({\bf x})| = \mathcal{O}(h_{1,T}^{k}), \label{11}\\
&  &  \sup_{{\bf x} \in \mathcal{C}}|\widehat {\widehat
m}_{\psi,T,add} ({\bf
x}) - \widehat {\widehat m}_{\psi,T,add} (t({\bf x}))| =\mathcal{O}\Big( \varepsilon(T)\Big)\label{12}\\
&  & \sup_{{\bf x} \in \mathcal{C}}|E\widehat {\widehat
m}_{\psi,T,add} ({\bf
x}) - E\widehat {\widehat m}_{\psi,T,add} (t({\bf x}))| \mathcal{O}\Big( \varepsilon(T)\Big),\label{13}\\
&  & \sup_{{\bf x} \in \mathcal{C}}|\widehat {\widehat
m}_{\psi,T,add} (t({\bf x})) - E\widehat {\widehat m}_{\psi,T,add}
(t({\bf x}))| = \mathcal{O}\Bigg(\Big(\frac{\log
T}{T}\Big)^{k/2k+1}\Bigg) ~~a.s., \label{14}\\
& & \sup_{{\bf x} \in \mathcal{C}}|\widehat {\widehat
m}_{\psi,T,add} ({\bf x}) - m_{\psi}({\bf x})|
=\mathcal{O}\Bigg(\Big(\frac{\log T}{T}\Big)^{k/(2k+1)}\Bigg)~~a.s..
\label{15}
\end{eqnarray}
\end{lemma}
{\it Proof of \ref{11}:} We have
\begin{eqnarray*}
&& \hspace{-2cm}E\widehat {\widehat m}_{\psi,T,add}({\bf x}) -
m_{\psi}({\bf x})
\\&=& \sum_{l=1}^d (E\widehat{\widehat \eta}_l(x_l) -\eta_l(x_l) )+
E(\int_{\mathbb{R}^d} \widetilde {\widetilde m}_{\psi,T}({\bf x})
q({\bf x}) d{\bf x}) - \int_{\mathbb{R}^{d}} m_{\psi}({\bf
x})q({\bf x}) d{\bf x}.
\end{eqnarray*}
By Fubini's theorem, we obtain
\begin{equation}\mbox{Bias}(\widehat{\widehat m}_{\psi,T,add}({\bf x}))=
\sum_{l=1}^d \mbox{Bias} (\widehat{\widehat \eta}_l(x_l)) +
\int_{\mathbb{R}^d} \mbox {Bias}(\widetilde{\widetilde
m}_{\psi,T}({\bf x})) q({\bf x})d {\bf x}. \label{decomp_bias1}
\end{equation}
We can write,
\begin{eqnarray}
\mbox{Bias}(\widehat{\widehat \eta}_1(x_1)) & = &
E(\widehat{\widehat \eta}_1(x_1)) - \eta_1(x_1)
\label{def_bias_eta}\nonumber \\
& = &  \{E( \widehat \alpha_1(x_1)) - \widetilde m_{\psi,T}(x_1)\}+E(\widehat C_T - C_T -C)\nonumber \\
&:= &(I) + (II).\label{decomp_bias_eta}
\end{eqnarray}
First consider the term $(I)$, we have
\begin{eqnarray*}
\widetilde m_{\psi,T}(x_1) &=&\frac{1}{T}\int_0^T E \Big\{
E(\psi(Y_t)|{\bf X}_t) \frac{\mathcal{G}({\bf
X}_t)}{f(X_{t,-1}|X_{t,1})} d{\bf x}_{-1}
\Big| X_{t,1}=x_1 \Big\}dt \\
&=& \int_{\mathbb{R}^{d-1}} m_\psi(x_1,{\bf
u}_{-1})\int_{\mathbb{R}^{d-1}} \frac{1}{h_{2,T}^{d-1}} K_{2}
\left(\frac{{\bf x}_{-1}- {\bf u}_{-1}}{h_{2,T}}\right)
q_{-1}({\bf
x}_{-1}) d{\bf x}_{-1}\ d{\bf u}_{-1}\\
&=&\int_{\mathbb{R}^{d-1}} m_\psi(x_1,{\bf u}_{-1})
\mathcal{G}({\bf u}_{-1}) d{\bf u}_{-1}.
\end{eqnarray*}
It follows that, under the conditions $(C.2)$, $(K.2)$ and $(K.3)$
\begin{eqnarray*}
&  &E(\widehat \alpha_1(x_1) ) - \widetilde m_{\psi,T}(x_1) =
\int_{\mathbb{R}}\frac{1}{h_{1,T}}
\widetilde m_{\psi,T}(u_1) K_1\Big( \frac{x_1-u_1}{h_{1,T}} \Big) du_1 -\widetilde m_{\psi,T}(x_1)\\
&=& \int_{\mathbb{R}} \big[\widetilde
m_{\psi,T}(x_1-v_1h_{1,T})-\widetilde m_{\psi,T}(x_1)\big]
K_1(v_1)dv_1\\
& =&\int_{\mathbb{R}} \int_{\mathbb{R}^{d-1}}
\Big[\sum_{i=1}^{k-1} \frac{(-h_{1,T}v_1)^i}{i!} \frac{\partial^i
m_{\psi}}{\partial x_1^i} (x_1,{\bf u}_{-1})\\
&& + \frac{(-h_{1,T}v_1)^k}{k!} \frac{\partial^k
m_{\psi}}{\partial x_1^k}(x_1-\theta h_{1,T} v1,{\bf u}_{-1}
)\Big]
\mathcal{G}({\bf u}_{-1}) d{\bf u}_{-1}  K_1(v_1)dv_1\\
& = & \int_{\mathbb{R}} \int_{\mathbb{R}^{d-1}}
\left[\frac{(-h_{1,T}v_1)^k}{k!} \frac{\partial^k
m_{\psi}}{\partial x_1^k}(x_1-\theta h_{1,T} v_1,{\bf u}_{-1}
)\right] \mathcal{G}({\bf u}_{-1}) d{\bf u}_{-1}  K_1(v_1)dv_1.
\end{eqnarray*}
Thus, we obtain
\begin{eqnarray}
&  &\sup_{x_1}|E(\widehat \alpha_1(x_1) ) - \widetilde
m_{\psi,T}(x_1)| = \mathcal{O}\left(h_{1,T}^k\right).
\label{res_decomp2_1}
\end{eqnarray}
Next, turning our attention to $(II)$, by (\ref{ii}) we have
\begin{eqnarray}
E(\widehat C_T -C_T) - C & =& \mathcal{O}(h_{1,T}^k).
\label{res_decomp2_2}
\end{eqnarray}
Combining (\ref{res_decomp2_1}) and (\ref{res_decomp2_2}), it
follows that
\begin{equation}
 \sup_{x_1}|E(\widehat{\widehat \eta}_1(x_1)) -
\eta_1(x_1)| = \mathcal{O}\left(h_{1,T}^k\right).
\label{res_bias_eta}
\end{equation}
On the other hand, we have, for all $0<\theta <1$,
\begin{eqnarray}
\mbox{Bias} (\widetilde{\widetilde m}_{\psi,T}({\bf x})) & := & E\widetilde{\widetilde m}_{\psi,T}({\bf x}) - m_\psi({\bf x}) \label{def_bias_mpsi}\\
& = & \int_{\mathbb{R}^d} [ m_\psi({\bf x} + h_{1,T}{\bf v})
-m_\psi({\bf x}) ]
K_3( {\bf v} )d{\bf v}.\nonumber\\
& = & \int_{\mathbb{R}^d} \sum_{i_1+...+i_d=k}
\frac{h_{1,T}^k}{k!}\frac{\partial^{i_1+...+i_d} m_\psi}{\partial
x_1^{i_1}...\partial x_d^{i_d}}({\bf x}+h_1\theta{\bf v})
v_1^{i_1}...v_d^{i_d} K_3({\bf v}) d{\bf
v}\\
&:=&\mathcal{O}(h_{1,T}^k), \label{res_bias_mpsi}\nonumber
\end{eqnarray}
Combining the decomposition (\ref{decomp_bias1}) and the statements
(\ref{res_bias_eta}) and (\ref{res_bias_mpsi}), we deduce the result
(\ref{11}). \vskip5pt \noindent{\it Proof of (\ref{12})}
 Under  the  condition $(K.5)$, there exists a constant M such that
\begin{eqnarray*}
\frac{1}{T}\int_0^T |Z_t({\bf x}) - Z_t({\bf t(x)})| dt\leq
\sum_{l=1}^d \frac{M}{h_{1,T}^2}|x_l - t({\bf x})_l|
\end{eqnarray*}
Consequently, using the expression of $r(T)$, we obtain
\begin{eqnarray*}
\sup_{{\bf x}\in \mathcal{C}}|\widehat{\widehat m}_{\psi,T,add}
({\bf x}) - \widehat {\widehat m}_{\psi,T,add} (t({\bf x}))|
=\mathcal{O} \big( \epsilon(T)\big).
\end{eqnarray*}\vskip5pt
\noindent{\it Proof of (\ref{13}):} Similarly as above, we may
deduce (\ref{13}).\vskip5pt \noindent {\it Proof of (\ref{14}):} In
view of Lemma \ref{lem3}, it is sufficient to prove discrete version
of (\ref{14}), that is
\begin{eqnarray}
\sup_{{\bf x} \in \mathcal{C}}|\widehat {\widehat m}_{\psi,n,add}
(t({\bf x})) - E\widehat {\widehat m}_{\psi,n,add} (t({\bf x}))| =
\mathcal{O}\Big( \varepsilon(n)\Big) ~~a.s.\label{discr}
\end{eqnarray}

Set  $n$ in $\mathbb{N}$ and,introduce some notations. Set,
\begin{eqnarray}
\widehat{\widehat m}_{\psi,n,add} ({\bf x}) - E\widehat{\widehat
m}_{\psi,n,add} ({\bf x}) &=:& \frac{1}{n} \int_0^n \xi_t({\bf
x})dt, \label{pour_12}
\end{eqnarray}
where
\begin{eqnarray*} \xi_t = \xi_t({\bf u}) := (Z_t({\bf u}) -
E(Z_t({\bf u})))
\end{eqnarray*}
and
\begin{eqnarray*}
Z_t & =& Z_t({\bf u}) \\
& =& \frac{\psi(Y_t)}{h_{1,n}h_{2,n}^{d-1}f({\bf X}_t)} \sum_{l=1}^d
\left\{ \int_{\mathbb{R}^{d-1}}
K_1\left(\frac{u_l-X_{t,l}}{h_{1,n}}\right)K_{2}\left(\frac{{\bf
x}_{-l}- {\bf
X}_{t,-l}}{h_{2,n}}\right)q_{-l}({\bf x}_{-l}) d{\bf x}_{-l}\right. \\
&  & -\left.\int_{\mathbb{R}^{d}}
K_1\left(\frac{x_l-X_{t,l}}{h_{1,n}}\right)K_{2}\left(\frac{{\bf
x}_{-l}- {\bf
X}_{t,-l}}{h_{2,n}}\right)q({\bf x}) d{\bf x}\right\} \\
&  & + \frac{1}{h_{1,n}^d}\frac{\psi(Y_t)}{f({\bf
X}_t)}\int_{\mathbb{R}^{d-1}} K_3\left(\frac{{\bf x}-{\bf
X}_t}{h_{1,n}}\right) q({\bf x}) d{\bf x}.
\end{eqnarray*}
Finally, we use the notation
\begin{eqnarray} V_i^n({\bf x}):=\frac{1}{n}\int_{(i-1)p}^{ip}\xi_t({\bf x})dt~~~i=1,...,2q' ~
\mbox{where}~p:=\frac{n}{2q'}:=\epsilon(n)^{-1/2}
 \label{p}\end{eqnarray}
So we can write
 \begin{eqnarray} \widehat{\widehat m}_{\psi,n,add}
({\bf x}) - E\widehat{\widehat m}_{\psi,n,add} ({\bf x}) & = &
\sum_{i=1}^{2q'}
  V_i^n({\bf x}).\label{pour_12}
\end{eqnarray}

We have to show that the following quantity is summable
\begin{eqnarray}
\mathbf{P}(\sup_{{\bf x} \in \mathcal{C}} |
\sum_{i=1}^{2q'}V_i(t({\bf x}))| \geq \varepsilon(n)) & \leq & r(n)
\sup_{p=1,...,r(n)} \mathbf{P}(| \sum_{i=1}^{2q'}V_i({\bf t}_p)|
\geq \varepsilon(n)). \label{inter1}
\end{eqnarray}
Let  $j$ be fixed  in $[1,r(n)]$. We have
$$
\mathbf{P}(| \sum_{i=1}^{2q'}V_i({\bf t}_j)| \geq
\varepsilon(n))\leq\mathbf{P}(| \sum_{i=1}^{q'}V_{2i}({\bf t}_j)|
\geq \varepsilon(n)/2)+\mathbf{P}(| \sum_{i=1}^{q'}V_{2i-1}({\bf
t}_j)| \geq \varepsilon(n)/2).
$$
Observing that for a given $M''$, $\xi_t({\bf
x})(\omega)<\frac{M''}{h_{1,n}},~\forall\omega\in\Omega$, we can use
recursively Bradley's lemma and define the independent random
variables $W_{2}({\bf t}_j),...,W_{2q'}({\bf t}_j)$ such that,
$\forall i\in[1,q']$, $W_{2i}$ and $V_{2i}$ have the same law and
$\forall \nu>0$
\begin{eqnarray}
P\bigg(|W_{2i}({\bf t}_j)-V_{2i}({\bf t}_j)|>\nu\bigg)&\leq&
11\bigg(\frac{\|V_{2i}({\bf
t}_j)\|_\infty}{\nu}\bigg)^{\frac{1}{2}}\alpha(p) \leq
11\bigg(\frac{pM''}{h_{1,n}\nu}\bigg)^{\frac{1}{2}}\alpha(p).
\label{bradley}
\end{eqnarray}
We have, for all $ 0<{\large\lambda} <
{\Large\frac{\varepsilon(n)}{2}}$
\begin{eqnarray}
\bigg\{|\sum_{i=1}^{q'}V_{2i}({\bf t}_j)|>\frac{\epsilon(n)}{2}\bigg\}&\subset&\bigg\{\big\{|\sum_{i=1}^{q'}V_{2i}({\bf t}_j)|>\frac{\epsilon(n)}{2};|V_i({\bf t}_j)-W_i({\bf t}_j)|\leq\frac{\lambda}{q'}~~~1\leq i\leq q'\big\}\bigg\}\nonumber\\
&&\bigcup\bigcup_{j=1}^{q'}\big\{|V_i({\bf t}_j)-W_i({\bf
t}_j)|>\frac{\lambda}{q'}\big\}\bigg\}. \nonumber
\end{eqnarray}
The choice {\large$\lambda$}={\Large$\frac{\epsilon(n)}{4}$} gives
us
\begin{eqnarray}
 P\bigg(|\sum_{i=1}^{q'}V_{2i}({\bf t}_j)|>\frac{\epsilon(n)}{2}\bigg)&\leq& P\bigg(
|\sum_{i=1}^{q'}W_{2i}({\bf t}_j)|>\frac{\epsilon(n)}{4}\bigg)
\nonumber\\&&+\sum_{i=1}^{q'}P\bigg(|V_{2i}({\bf t}_j)-W_{2i}({\bf
t}_j)|>\frac{\epsilon(n)}{4q'}\bigg).\label{sep}
\end{eqnarray}
We treat separately the two terms of the last inequality. For the
second one, the application of (\ref{bradley})under the condition
$(A.1)$ drives us to
\begin{eqnarray}
\sum_{i=1}^{q'}P\bigg(|V_{2i}({\bf t}_j)-W_{2i}({\bf
t}_j)|>\frac{\epsilon(n)}{4q'}\bigg) \leq 11q'
\Big(\frac{4q'pM'}{h_{1,n}\epsilon
(n)}\Big)^{1/2}\alpha(p)=\mathcal{O}(r(n)^{-1}n^{\mu})
\label{borela} \\
~where~\mu<-1.\nonumber\end{eqnarray} In order to dominate
$P\bigg( |\sum_{i=1}^{q'}W_{2i}({\bf
t}_j)|>\frac{\epsilon(n)}{4}\bigg)$, we must bound the variance of
$W_{2i}$ (which has the same law as $V_{2i}$) to use Bernstein's
inequality
\begin{eqnarray}
{\rm Var}(W_{2i}({\bf t}_j))&=& E\big(V_{2i}({\bf t}_j)^2\big)\nonumber \\
&\leq&\frac{1}{n^2}\int_{[(2j-1)p,2jp]}E(\xi_t^2)dt
\end{eqnarray}
The kernels are bounded, so we can easily see, after a change of
variables, that there exists a constant $M'''$ such that
$$E(Z_t^2)\leq\frac{M'''}{h_{1,n}},$$
witch implies
$$\mathbb{E}\big(\xi_t^2\big)\leq\frac{M'''}{h_{1,n}}~\mbox{and}~{\rm Var}(W_{2i}({\bf t}_j))\leq\frac{pM'''}{n^2h_{1,n}}.$$
Observe that, for a given $S$ in $\mathbb{R}^{*+}$,
$\xi_t(\omega)<\frac{S}{h_{1,n}},~\forall\omega\in\Omega$, we
readily have $$E|W_i|^k\leq (p\frac{M'}{nh_{1,n}})^{k-2}p!E|W_i|^2,
\forall i.$$ This allows us to apply Bernstein's inequality
\begin{eqnarray} P\bigg(|\sum_{i=1}^{q'}W_i({\bf t}_j)|>
\frac{\epsilon(n)}{4}\bigg)&\leq&
2\exp\bigg(-\frac{\epsilon(n)^2}{16(\frac{4q'pM}{n^2h_{1,n}}+\frac{M'p\epsilon(n)}{2nh_{1,n}})}\bigg)
\nonumber\\&=& 2
\exp\bigg(-\frac{\epsilon(n)^2nh_{1,n}}{32M+8M'p\epsilon(n))}
\bigg). ~~~~~~~~~~~~~~~~~~~~~~~\label{prob}\end{eqnarray} The
expression of $p$ and $\varepsilon(n)$ gives us
$p\varepsilon(n)\rightarrow0$ and the sequence $\sum_{n=1}^N
r(n)P\bigg(|\sum_{i=1}^{q'}W_i({\bf t}_j)|>
\frac{\epsilon(n)}{4}\bigg)$ converges as N grows to infinity if
we choose a large enough $C$ in $\epsilon(n)$. In view of this last inequality and (\ref{borela}), we obtain (\ref{discr}) by Borel-Cantelli.\\
\\
{\it Proof of (\ref{15}):}~~ By (\ref{estcarbo}), we have
\begin{eqnarray*}
\sup_{{\bf x}\in \mathcal{C}} \big|\widetilde {\widetilde
m}_{\psi,T}({\bf x})- {\widetilde m}_{\psi,T}({\bf x})\big|
&\leq&M \frac{ \sup_{{\bf x}\in \mathcal{C}} \big|f({\bf x})-
\widehat f_T( {\bf x})\big|}{\inf_{{\bf x}\in \mathcal{C}}
f^2({\bf x}) +o(1)}
\frac{1}{Th_{1,T}^d}\int_0^T\Big|K_{3}\Big({\frac{{\bf x} -{\bf
X}_{t}}{h_{1,T}}} \Big)\Big| dt.
\end{eqnarray*}
Using  the statements (\ref{additive_component}) and
(\ref{additive_component_marginale}), and the Theorem on a density
estimator due to \cite{Bosq1993}, we obtain
\begin{eqnarray*} &&\hspace{-2cm} \sup_{{\bf x} \in \mathcal{C}}
|\widehat{\widehat m}_{\psi,T,add}({\bf x}) - {\widehat
m}_{\psi,T,add}({\bf x})|\\ &\leq & 2d\max_{1 \leq l \leq d}
\sup_{{\bf x}\in \mathcal{C}} \big|\widetilde {\widetilde
m}_{\psi,T,l}({\bf x})- {\widetilde m}_{\psi,T,l}({\bf x})\big| +
\sup_{{\bf x}\in \mathcal{C}} \big|\widetilde {\widetilde
m}_{\psi,T}({\bf x})- {\widetilde m}_{\psi,T}({\bf x})\big|\\
& =&\mathcal{O}\left(\Big(\frac{\log
T}{T}\Big)^{k/(2k+1)}\right)~~\mbox{a.s.}.
\end{eqnarray*}

%

\end{document}